\documentclass[12pt]{article}
\usepackage{amssymb}
\usepackage{amsthm}

\newtheorem{Theorem}{Theorem}
\newtheorem{Proposition}{Proposition}
\newtheorem{Lemma}{Lemma}
\newtheorem{Question}{Question}
\newtheorem{Conjecture}{Conjecture}

\begin{document}
\title{A few remarks about the Hilbert scheme \\ of smooth projective curves}
\author{Edoardo Ballico and Claudio Fontanari}
\date{}

\maketitle

\begin{abstract} We discuss two conjectures by Francesco Severi and Joe Harris about the irreducibility and the dimension of the Hilbert scheme parameterizing smooth projective curves of given degree and genus. 

\vspace{0.3cm}



\end{abstract}

\section{Introduction}

\emph{Man erkennt zun\"achst \emph{(...)}, da{\ss} die irreduziblen Kurven $C$ vom Geschlecht $p$ 
und der Ordnung $n \ge p+r$ des Raumes $S_r$ eine einzige Familie $V$ bilden, deren allgemeine Kurve 
nicht spezial ist}. This statement made by Francesco Severi at the beginning of \S 5 of Anhang G in 
\cite{V} (see p. 368), turns out to be false, as an example suggested by Joe Harris first showed 
in the mid eighties (see \cite{E2}, Proposition 9). Here we focus on the problem  
of establishing the range in which Severi's claim holds true.   

Let $\mathcal{H}_{d,g,n}$ denote the open subset of the Hilbert scheme parameterizing 
smooth, irreducible, nondegenerate curves of degree $d$ and genus $g$ in $\mathbb{P}^n$.   

\begin{Theorem}\label{main}
Fix integers $g \ge 0$, $n \ge 2$ and $d \ge g+n$. If $5 \le n \le 6$ assume $d > \frac{2n+1}{n+5} g + 1$.  
If $n \ge 7$ assume $d > \frac{2n-4}{n+1} g + \frac{n+13}{n+1}$. Then $\mathcal{H}_{d,g,n}$ is irreducible. 
\end{Theorem}

The above statement collects several already known results together with a few original ones. 
The case $n=2$ was first pointed out by Arbarello and Cornalba (see \cite{AC}, Lemma 3.4), 
while the cases $n=3,4$ are due to Ein (see \cite{E1}, Theorem 4, and \cite{E2}, Theorem 7; 
see also \cite{KK2}, Theorem 1.5, for a different proof of the case $n=3$ in the spirit of 
\cite{AC}). The result for $n=5$ is due to Iliev (see \cite{I2}, Theorem A), while the  
statement for $n=8, g \ge 19$ and for $n \ge 9$ has been proven by Kim (see \cite{Kim}, 
Theorem 3.7). To the best of our knowledge, the cases $n = 6$, $n = 7$, and $n = 8, 
g \le 18$ have not previously appeared in the literature (the sharpest published bound 
covering this range we are aware of is indeed $d > \frac{2n-2}{n+2} g + \frac{n+8}{n+2}$, 
see \cite{E2}, Theorem 8). Our proof extends to every $n$ the lines of the argument 
provided by Iliev in \cite{I2} for the case $n=5$. 

The hypothesis $d > \frac{2n-4}{n+1} g + \frac{n+13}{n+1}$ in the statement of Theorem \ref{main}
turns out to be sharp (see for instance \cite{MS}). On the other hand, the stronger numerical 
assumption $d > \frac{2n+1}{n+5} g + 1$ for $5 \le n \le 6$ seems to be just a technical one, 
even though we were not able to remove it. Our partial results in this direction are 
summarized in the following:

\begin{Proposition}\label{nonvanishing}
Fix integers $g \ge 0$, $n \ge 2$ and $d \ge g+n$. If $d > \frac{2n-4}{n+1} g + \frac{n+13}{n+1}$
and $\mathcal{H}_{d,g,n}$ is not irreducible, then $H^0(C, K_C - N_{L-D} - D) \ne 0$, where 
$C$ is a general curve in any extra component of $\mathcal{H}_{d,g,n}$, $L = \mathcal{O}_C(1)$, 
$h^0(C,L)=r+1$, $D = p_1 + \ldots + p_{r-3}$ with $p_i$ general points on $C$, and $N_{L-D}$ 
denotes the normal bundle to $C$ in the immersion into $\mathbb{P}^3$ defined by the complete 
linear series associated to $L-D$.    
\end{Proposition}

By applying \cite{E2}, Theorem 5, to $N_{L-D}$ we obtain $h^0(C, K_C - N_{L-D}) 
\le g-d+r-2$. If $D$ were general with respect to $L-D$, then we could conclude 
that $h^0(C, K_C - N_{L-D} - D) = 0$, but unluckily this is not the case. 

We also investigate the range in which the Hilbert scheme is reducible, in particular 
in Section \ref{harris} we provide the following counterexample to a conjecture stated 
in \cite{H} (see also \cite{EH}, question (1) on p. 40): 

\begin{Theorem}\label{counterexample}  
For every integer $g \ge 12$, let $g \equiv \varepsilon$ modulo $4$ with $0 \le \varepsilon \le 3$. 
Then there exist integers $d, n$ and a component $\mathcal{H}$ of the Hilbert scheme whose general 
member corresponds to a smooth, irreducible, nondegenerate curve of degree $d$ and genus $g$ 
in $\mathbb{P}^n$ such that $\dim \mathcal{H} > 3 g - 3 + \rho(g, n, d) + (n+1)^2-1$ and the image 
of the rational map $\mathcal{H} \to \mathcal{M}_g$ has codimension $\frac{g}{2} - \frac{\varepsilon}{2} + 2$. 
\end{Theorem}

As pointed out in \cite{H}, a more refined version of the same conjecture predicts the existence of 
a number $\beta(g)$ tending linearly to $\infty$ with $g$, such that any such component $\mathcal{H}$ 
whose image in $\mathcal{M}_g$ has codimension $\beta \le \beta(g)$ has the expected dimension. 
From this point of view, Theorem \ref{counterexample} only implies that $\beta(g) \le \frac{g}{2} + 2$, 
and we believe that such a problem is well worth of further investigation. 

This research project was partially supported by MIUR and GNSAGA of INdAM. The second named author 
started thinking about the Hilbert scheme in the warm atmosphere of the 2011 Summer Session of 
IAS/Park City Mathematical Institute, whose support is gratefully acknowledged. In particular, 
he benefitted from encouraging discussions with Joe Harris on Brill-Noether theory. His interest  
in the irreducibility problem was aroused by the beautiful survey \cite{S} of Edoardo Sernesi, 
to whom he is also indebt for stimulating pleasant conversations about Severi's Anhang G.  

We work over the complex field. 

\section{Irreducibility}

Fix nonnegative integers $g, r, d$ and let $\mathcal{W}^r_d = \{ (C,L): C \textrm{ a smooth 
curve of}$ $\textrm{genus $g$, $L$ a line bundle on $C$ of degree $d$ with } h^0(C,L)=r+1 \}$
(for a rigorous definition we refer to \cite{ACG}, Chapter XXI., \S 3). 

We borrow from \cite{I2} the three main technical ingredients of the proof of Theorem \ref{main}. 

\begin{Lemma} \label{dimension}
Assume $r \ge 2$ and $g - d + r \ge 2$ and let $\mathcal{W}$ be an irreducible component 
of $\mathcal{W}^r_d$ with general element $(C,L)$. If the moving part of $\vert L \vert$ is

(i) very ample, then $\dim \mathcal{W} \le 3d + g + 1 - 5r$. 

(ii) birationally very ample, then $\dim \mathcal{W} \le 3d + g - 1 - 4r$. 

(iii) compounded, then $\dim \mathcal{W} \le 2g - 1 + d - 2r$.   
\end{Lemma}

\proof For (i), see \cite{E2}, Theorem 6. For (ii), see \cite{KK1}, Lemma 1. For (iii), see \cite{I1}, 
Proposition 2.1. 

\qed

\begin{Lemma} \label{strictly}
Assume $r \ge 2$ and $g - d + r \ge 3$ and let $\mathcal{W}$ be an irreducible component 
of $\mathcal{W}^r_d$ whose general element $(C,L)$ is such that $L$ is strictly birationally 
very ample and the moving part of $K_C - L$ is birationally very ample. Then there exists 
an irreducible component $\mathcal{W}'$ of $\mathcal{W}^{r-1}_{d-2}$ such that 
$\dim \mathcal{W} \le \dim \mathcal{W}'$ and if $(C',L')$ is the general element 
of $\mathcal{W}'$ then the moving part of $K_{C'} - L'$ is birationally very ample. 
\end{Lemma}

\proof See \cite{I2}, Lemma 2.2 and Remark 2.3. 

\qed

\begin{Lemma} \label{clifford}
Let $g^r_d$ be a birationally very ample linear series of degree $d \ge g$ on a smooth curve of genus $g$. 
Then 
$$
r \le \frac{1}{3} (2d-g+1)
$$
\end{Lemma}

\proof See either \cite{ACGH}, Chapter III., \S 2, (2.2) on p. 115, or \cite{E1}, Lemma 7, or \cite{I1}, Lemma 2.3.  

\qed

\emph{Proof of Theorem \ref{main}.} By \cite{E1}, Theorem 4, and \cite{E2}, Theorem 7, we may assume
\begin{equation}\label{five}
n \ge 5. 
\end{equation}
Since $d \ge g+n$, there is a unique open subset of $\mathcal{H}_{d,g,n}$ corresponding to nonspecial curves. 
Suppose by contradiction that $\mathcal{H}_{d,g,n}$ is reducible. Then there exists an irreducible component 
$\mathcal{H}$ of $\mathcal{H}_{d,g,n}$ such that the general curve $C$ in $\mathcal{H}$ satisfies 
$h^0(C, \mathcal{O}_C(1)) = r+1$ and $h^1(C, \mathcal{O}_C(1)) = \delta > 0$. In particular, if 
$\mathcal{Y}$ is the irreducible component of $\mathcal{W}^r_d$ such that $(C, \mathcal{O}_C(1)) \in 
\mathcal{Y}$, then we must have
\begin{equation}\label{existence} 
\dim \mathcal{Y} + (n+1)(r-n) \ge 4g-3-(n+1)(g-d+n). 
\end{equation}
First of all, we claim that 
\begin{equation}\label{three}
\delta \ge 3.
\end{equation}
Indeed, if $\delta \le 2$ then by \cite{E2}, Theorem 6 (a), $\dim \mathcal{Y} = 4g-3-(r+1)(g-d+r)$, 
hence from (\ref{existence}) it follows that $n=r$ and $0 < \delta =g-d+r=g-d+n \le 0$, contradiction.  

The idea now is to focus on $K_C - \mathcal{O}_C(1)$. Let $B$ be the base locus of $K_C - \mathcal{O}_C(1)$ 
with $\deg B = b$ and let $\mathcal{W}$ be the irreducible component of $\mathcal{W}^{\delta - 1}_{2g-2-d}$ 
such that $(C, K_C - \mathcal{O}_C(1)) \in \mathcal{W}$. 

Assume first that $K_C - \mathcal{O}_C(1) - B$ is very ample. In this case, Lemma \ref{dimension} (i) 
applied to $\mathcal{W}$ yields $\dim \mathcal{Y} = \dim \mathcal{W} \le 3(2g-2-d) + g + 1 - 5(\delta-1)$
and (\ref{existence}) implies 
\begin{equation}\label{bound}
(n-1)d \le (n-1)g + (n-4)r + 3.
\end{equation}
Hence by Lemma \ref{clifford} we deduce $d \le \frac{2n+1}{n+5} g + 1$, which contradicts our 
numerical assumptions. 

Assume next that $K_C - \mathcal{O}_C(1) - B$ is composed with an involution of degree $m$ onto a curve 
$\Gamma$ of genus $\gamma \ge 0$. If $f$ is the induced morphism, then $f = h \circ g$, with 
$g: C \to \Gamma$ of degree $m$ and $h: \Gamma \to \mathbb{P}^{\delta - 1}$ of degree
$\frac{2g-2-d-b}{m}$. It follows that $\frac{2g-2-d-b}{m} \ge \delta - 1$ and 
\begin{equation}\label{degree}
m \le \frac{2g-2-d-b}{\delta - 1} \le  \frac{2g-2-d}{\delta - 1}. 
\end{equation}

Let now $k := \frac{n-2}{3(n-1)}$. We claim that if $d > \frac{2n-4}{n+1} g + \frac{n+13}{n+1}$ 
then $(\delta - 1) > k(2g-2-d)$. Indeed, by Lemma \ref{dimension} (iii) applied to $\mathcal{W}$ 
we have $\dim \mathcal{Y} = \dim \mathcal{W} \le 2g-1+(2g-2-d)-2(\delta-1)$ 
and from (\ref{existence}) we deduce that $2g-n-3 \le (2g-2-d)+(n-1)(\delta-1)$.  
If $(\delta - 1) \le k(2g-2-d)$, then $2g-n-3 \le (1+(n-1)k)(2g-2-d)$ and since 
$d > \frac{2n-4}{n+1} g + \frac{n+13}{n+1}$ we reach a numerical contradiction. 

Hence (\ref{degree}) implies $m < \frac{1}{k} \le 4$, where the last 
estimate follows from (\ref{five}), so that $m \le 3$.  

If $m = 3$, then from (\ref{degree}) we also deduce
\begin{eqnarray}
b &\le& -3(\delta - 1)+(2g-2-d) \label{b}\\
(\delta-1) &\le& \frac{2g-2-d}{3} \label{delta}.
\end{eqnarray}

If $m = 2$, then $\gamma \ge 1$ since otherwise $C$ would be hyperelliptic 
and the dual of $K_C - \mathcal{O}_C(1)$ could not be very ample.
If $\frac{2g-2-d-b}{2} \le 2 \gamma$ then by Clifford's Theorem we obtain 
$(\delta -1) \le \frac{2g-2-d-b}{4}$, hence both (\ref{b}) and (\ref{delta}) 
hold.  
If instead $\frac{2g-2-d-b}{2} > 2 \gamma$ then the linear series 
$\vert H \vert$ of degree $\frac{2g-2-d-b}{2}$ on $\Gamma$ defining 
$h: \Gamma \to \mathbb{P}^{\delta - 1}$ is nonspecial. It follows 
that $\dim \vert H + Q \vert = \dim \vert H \vert + 1$ for every $Q \in g(C)$ 
and $\dim \vert K_C - \mathcal{O}_C(1) + P_1 + P_2 \vert \ge \dim 
\vert K_C - \mathcal{O}_C(1) \vert + 1$ for $P_1, P_2 \in C$ such that 
$g(P_1)=g(P_2)=Q$. Hence by Riemann-Roch we have 
$\dim \vert \mathcal{O}_C(1) - P_1 - P_2 \vert \ge 
\dim \vert K_C - \mathcal{O}_C(1) \vert + d - g = 
\dim \vert \mathcal{O}_C(1) \vert - 1$, contradicting 
the very ampleness of $\mathcal{O}_C(1)$.   

Now, for $\gamma = 0$ we have $\dim \mathcal{W} \le \dim \mathcal{W}^1_m + b = 2g+2m-5+b \le 2g+1+b$, 
while for $\gamma \ge 1$ we have $\dim \mathcal{W} \le 2g-2+(2n-3)(1-\gamma)+\gamma+b \le 2g+1+b$, 
hence in both cases by (\ref{b}) we get $\dim \mathcal{Y} = \dim \mathcal{W} \le 2g+1-3(\delta - 1)+(2g-2-d)$. 
From (\ref{existence}) we deduce $2g-5-n \le (2g-2-d)+(n-2)(\delta - 1)$ and by (\ref{delta}) 
we conclude $d \le \frac{2n-4}{n+1} g + \frac{n+13}{n+1}$, contradicting our numerical 
assumptions. 

Assume finally that $K_C - \mathcal{O}_C(1) - B$ is strictly birationally very ample. 
By an iterative application of Lemma \ref{strictly}, either we are reduced to one 
of the previous two cases, or we obtain a component $\mathcal{W}'$ of 
$\mathcal{W}^2_{2g-2-d-2(\delta-3)}$ such that the general element of 
$\mathcal{W}'$ is strictly birationally very ample and 
$\dim \mathcal{W} \le \mathcal{W}'$. Hence Lemma \ref{dimension} (ii) 
applied to $\mathcal{W}'$ yields $\dim \mathcal{Y} = \dim \mathcal{W} \le \dim \mathcal{W}' 
\le 3(2g-2-d-2(\delta-3))+g-9$ and from (\ref{existence}) it follows that 
\begin{equation}\label{one}
(n-1)d \le (n-1)g+(n-4)r - \delta + 6.
\end{equation}
On the other hand, since we know that (\ref{bound}) is incompatible with our numerical assumptions, we have 
\begin{equation}\label{two}
(n-1)d > (n-1)g+(n-4)r + 3.
\end{equation}
From (\ref{one}) and (\ref{two}) we deduce $\delta < 3$, which contradicts (\ref{three}). 

\qed 

\emph{Proof of Proposition \ref{nonvanishing}.} 
Assume by contradiction $H^0(C, K_C - N_{L-D} - D) \ne 0$. 
From the exact sequence (see \cite{BEL} (2.7) and (2.5))
$$
0 \to K - N_{L-D} - D \to K - N_L \to \sum_{i=1}^{r-3} K-L(-2p_i) \to 0
$$
we deduce 
$$
h^1(C, N_L) = h^0(C, K - N_L) \le (r-3)(\delta-2). 
$$
If $\mathcal{Y}$ is the irreducible component of $\mathcal{W}^r_d$ such that 
$(C, L) \in \mathcal{Y}$ then $\dim \mathcal{Y} \le 4g-3+(r+1)(g-d+r)+h^1(C, N_L)
\le 4g-3+(r+1)(g-d+r)+(r-3)(\delta-2)$ and from (\ref{existence}) it follows that 
$(n-3)d \le (n-3)g + (n-5)r + 6$. Hence by Lemma \ref{clifford} we deduce 
$d \le \frac{2n-4}{n+1} g + \frac{n+13}{n+1}$, 
which contradicts our numerical assumption.

\qed

\section{Dimension}\label{harris}

Fix nonnegative integers $g,n,d$ and let $\rho(g,n,d) = g-(n+1)(g-d+n)$ be the so-called Brill-Noether number.

In the recent survey \cite{H}, p. 142, the following conjecture is stated 
(see also \cite{EH}, \S 2.b., pp. 39--40): 

\begin{Conjecture}
If $\mathcal{H}$ is any component of the Hilbert scheme whose general member corresponds
to a smooth, irreducible, nondegenerate curve of degree $d$ and genus $g$ in 
$\mathbb{P}^n$, and the image of the rational map $\mathcal{H} \to \mathcal{M}_g$ has 
codimension $g - 4$ or less, then $\dim \mathcal{H} = 3 g - 3 + \rho(g, n, d) + (n+1)^2-1$.
\end{Conjecture}

Immediately after that (see \cite{H}, p. 143), it is remarked that \emph{to be 
honest, the available evidence suggests simply the existence of a number $\beta(g)$
tending linearly to $\infty$ with $g$, such that any such component $\mathcal{H}$ whose 
image in $\mathcal{M}_g$ has codimension $\beta \le \beta(g)$ has the expected 
dimension; we use the function $g - 4$ just for simplicity}.

Theorem \ref{counterexample} shows that the literal statement of the above conjecture 
turns out to be false and that the number $\beta(g)$ satisfies the inequality
$\beta(g) \le \frac{g}{2} + 2$. 

\emph{Proof of Theorem \ref{counterexample}.} 
Let $k = \frac{g+\varepsilon}{4}$ and let $C$ be a general $k$-gonal curve 
with $g^1_k = \vert E \vert$. By \cite{K}, Corollary 3.3, 
$\vert K - 2 E \vert$ is very ample, so the general $k$-gonal curve is embedded 
as a curve of degree $d$ in $\mathbb{P}^n$ with $d = 2g-2-2k$ and $n = g-2k+1$. 
Since $k < \frac{g}{4}+1$, we have $2g+2k-5 > 3g-3+\rho(g,n,d)$, hence the family 
of general $k$-gonal curves embedded by $\vert K - 2 E \vert$ belongs to a component 
$\mathcal{H}$ of the Hilbert scheme such that $\dim \mathcal{H} > 3g-3+\rho(g,n,d) + (n+1)^2-1$. 
On the other hand, the codimension in $\mathcal{M}_g$ of the family of 
general $k$-gonal curves is $3g-3-(2g+2k-5) = \frac{g}{2} - \frac{\varepsilon}{2}+2$.

\qed

We are wondering whether Theorem \ref{counterexample} could be improved or not:

\begin{Question}
At least for $g$ large enough, is $\beta(g) = \frac{g}{2} + 2$ and is the family 
above the only one of smallest possible codimension in $\mathcal{M}_g$?
\end{Question}

\vspace{0.3cm}

\noindent
Edoardo Ballico \newline
Dipartimento di Matematica \newline 
Universit\`a di Trento \newline 
Via Sommarive 14 \newline 
38123 Trento, Italy. \newline
E-mail address: ballico@science.unitn.it

\vspace{0.3cm}

\noindent
Claudio Fontanari \newline
Dipartimento di Matematica \newline 
Universit\`a di Trento \newline 
Via Sommarive 14 \newline 
38123 Trento, Italy. \newline
E-mail address: fontanar@science.unitn.it

\end{document}